# *Stable Outcomes For Contract Choice Problems*

by
**Somdeb Lahiri**
**School of Economic and Business Sciences,**
**University of Witwatersrand at Johannesburg,**
**Private Bag 3, WITS 2050,**
**South Africa.**
**November 2003.**
Email: lahiris@sebs.wits.ac.za
Or
lahiri@webmail.co.za

*Abstract*
In this paper, we consider the problem of choosing a set of multi-party contracts, where each coalition of agents has a non-empty finite set of feasible contracts to choose from. We call such problems, contract choice problems. We provide conditions under which a contract choice problem has a non-empty set of "stable" outcomes. There are two types of stability concepts we study in this paper: cooperative stability and non-cooperative stability. The cooperative stability concept that we invoke here is the core. The non-cooperative stability concept that we study here is individual stability. The final result of this paper states that every contract choice problem has a non-empty weak bargaining set.

1. Introduction:

In this paper, we consider the problem of choosing a set of multi-party contracts, where each coalition of agents has a non-empty finite set of feasible contracts to choose from. We call such problems, contract choice problems. The economic motivation behind the problem, arises from several real world "commons problems", where agents can pool their initial resources and produce a marketable surplus, which needs to be shared among themselves. There are clearly, two distinct problems that arise out of such real world possibilities: (i) Coalition Formation: Which are the disjoint coalitions that will form in order to pool in their resources? (ii) Distribution: How will a coalition distribute the surplus within itself? While, the possibility of an aggregate amount of surplus being generated by a coalition is fairly common, there are many situations where more than one aggregate surplus results from a cooperative activity, and the distribution of the surplus depends on the particular aggregate that a coalition chooses to share.
The model we propose is a generalization of the model due to Shapley and Scarf (1974) called the housing market. Shapley and Scarf (1974), considers a private ownership economy, where each individual owns exactly one object and what is sought is the existence of an allocation in the core of the economy. In our model each non-empty subset of agents has a non-empty finite set of pay-off vectors to choose from. An outcome comprises a partition of the set of agents, and an assignment for each coalition in the partition a feasible pay-off vector. Our model is therefore a special kind of cooperative game with non-transferable utility. In the context of our contract choice model, the Shapley-Scarf housing market corresponds to a situation, where each individual assigns a monetary worth to each object, and a feasible pay-off vector for a coalition, is the set of utility vectors available to the coalition, when it re-allocates objects within itself, without any one in the coalition retaining his initial endowment, unless the coalition is a singleton.

Our primary objective in this paper, is to provide conditions under which a contract choice problem has a non-empty set of "stable" outcomes. There are two types of stability concepts we study in this paper: cooperative stability and non-cooperative stability. The cooperative stability concept that we invoke here is the core. An outcome is said to belong to the core of a contract choice problem, if there is no subset of agents who could select a feasible pay-off and be better off. Roth and Postelwaite (1977) used Gale's *Top Trading Cycle Algorithm* to show that if preferences are strict, then there exists a unique competitive equilibrium allocation, which is also the unique core allocation, for a Shapley-Scarf housing market. However, we are able to show with the help of a three agent example, that there exists a generalized contract choice problem, which does not admit any stable outcome.

We show here, that an adaptation of the weak top coalition property due to Banerjee, Konishi and Sonmez (2001), guarantees the non-emptiness of the core. The original version of the property due to Banerjee, Konishi and Sonmez (2001) that we adapt here, were postulated for coalition formation games, and as such do not apply to our context of a contract choice problem.

A salient feature of many markets is to match one agent with another. This is particularly true, in the case of assigning tasks to individuals where each task is under the supervision of an individual, and where the set of supervisors and the set of workers are disjoint. Such markets are usually studied with the help of "two sided matching models" introduced by Gale and Shapley (1962) called the marriage problem. The solution concept proposed by Gale and Shapley (1962), called a stable matching, requires that there should not exist two agents, who prefer each other, to the individual they have been paired with. It was shown in Gale and Shapley (1962), in a framework where every agent has preference defined by a linear order over the entire set of agents, that a marriage problem always admits a stable matching. An overview of the considerable literature on marriage problems that has evolved out of the work of Gale and Shapley (1962), is available in Roth and Sotomayor (1990). Lahiri (2002) contains alternative simpler proofs of some existing results and some new conclusions for two-sided matching problems.

In Lahiri (2003 b), we propose a generalization of a model due to Eriksson and Karlander (1998). We allow each pair of agents a non-empty finite set of real valued divisions of a good to choose from. Each agent is assumed to prefer more of the good to less of it. Further, the set of agents are divided into two disjoint sets, with one set being the set of men and the other the set of women, with no pair of the same sex being able to obtain an allocation which is at least as good as an allocation that could be obtained by them remaining single or by forming a pair with a member of the opposite sex. If each pair of agents is provided singletons to choose from, then we have the marriage problem of Gale and Shapley (1962). In Lahiri (2003 c) we show that the generalization proposed in Lahiri (2003 b) is simply a particular type of contract choice problem, which invariably admits a non-empty core. We show, that a simple generalization of the Deferred Acceptance Procedure with men proposing due to Gale and Shapley (1962), yields outcomes for the generalized marriage problem, which necessarily belong the core. We also show, that any outcome of this procedure is Weakly Pareto Optimal for Men, i.e. there is no other outcome which all men prefer to an outcome of this procedure. This result is an extension to our framework, of a similar result due to Roth and Sotomayor (1990).

The non-cooperative stability concept that we study here is individual stability, which has been proposed by Bogomolnaia and Jackson (2002) for hedonic coalition formation games, "where no allocation of goods need to be kept track of". Our individual stability concept is a modification of the one due to Bogomolnaia and Jackson (2001), defined so that it is consistent with the framework of our analysis. Bogomolnaia and Jackson (2001), motivate the applicability of their individual stability concept by providing the examples of "professors changing universities, soccer players considering changing teams,…, individuals changing clubs". However, there may be professors who select universities, on the basis of the net remuneration package that each university has to offer them, rather than on the identity of the individuals they would be associated with were they to be employed by a particular university. Similarly, a soccer player may choose a club, only if that club provides him at least the same remuneration that he was receiving in his present assignment. In such situations the value to an individual of belonging to a coalition is determined by what the coalition assigns to the individual and not merely by the identity of other members of the coalition.

We say that a feasible outcome is individually stable, if there is no agent who can unilaterally deviate by joining a coalition of agents who were earlier part of a group, and thereby improve the condition of every member of the coalition he joins (including himself!). Clearly an outcome in the core is individually stable, although the converse need not be true. We show that a property, referred to here as *weak top cycle property*, suffices to guarantee the existence of an individually stable allocation. The weak top cycle property says that given any non-empty subset of agents V, there exists a non-empty subset S of V containing s distinct agents, an outcome and a one to one function $\psi$ from S to the set $\{1,\ldots,s\}$ such that: (a) members of S form a coalition at the outcome and receive a pay-off vector at the outcome that is feasible for S; (b) each agent in S, prefers the stated outcome to any that he would be getting by forming a coalition with one or more agents in V who do not belong to S; (c) given any two agents a and b in S, if $\psi(a) < \psi(b)$, and for some coalition contained in S there is a pay-off vector, which 'a' prefers to the given outcome, then some other member of the coalition prefers the outcome to the pay-off vector. It is worth emphasizing that the weak top cycle property, has some resemblance to the concepts of consecutive NTU games due to Greenberg and Weber (1986) and consecutive coalition formation games due to Bogomolnaia and Jackson (2002).

Zhou (1994) introduced a concept of the bargaining set, which is a slight variation of the original one due to Aumann and Maschler (1964). Yet another notion of a bargaining set is due to Mas-Colell (1989). The Zhou(1994) bargaining set of a marriage problem always contains its non-empty core. Klijn and Masso( undated) introduced the concept of the weakly stable set for a marriage problem and showed that it coincided with its bargaining set as defined by Zhou (1994).

In a final section of this paper we introduce the concepts of the weak bargaining set for contract choice problems. Our concepts resemble a possible extension of similar concepts for marriage problems, due Klijn and Masso (undated). The basic idea behind the weak bargaining set is a set of feasible allocations which do not admit a credible objection (i.e. every strong objection has a strong counter-objection). Our definition of a credible objection is somewhat different from that of Zhou (1994) or Mas-Colell (1989), in that we require a strong counter-objection to make none of its proponents worse off than what they were at the time when the objection was raised. We further require that no sub-coalition of an objecting coalition can block the objecting pay-off. We show by a three agent example, that a natural analog of the bargaining set due to Mas-Colell (1989) may well be empty for room-mates problems. The final result of this paper states that every contract choice problem has a non-empty weak bargaining set.

A related paper [Lahiri (2003 a)] studies conditions which guarantee the existence of 'stable' allocations in a generalized matching model, where each of a finite number of agents owns a single indivisible objects, which can be re-allocated among them, so long as the resulting allocation is not *a priori* infeasible. Each agent has a strict ranking over the set of indivisible objects. In that paper generalized matching problems are referred to as a housing market, and sufficient conditions for the existence of non-empty cores and individually stable sets provided, which bear some resemblance to ones provided here. The marriage and room-mates problems of Gale and Shapley (1962) are special cases of this model.

2. Contract Choice Problems: Let X be a non-empty finite subset of $\aleph$ (: the set of natural numbers), denoting the set of participating agents. We assume that each agent prefers more money to less.

**Given $S \in [X]$, let $C(S) = \{\mu / \mu$ is a bijection on X with $\mu(S) = S\}$ and $C^0(S) = \{\mu \in C(S)$ / T is a non-empty proper subset of S implies $\mu(T) \neq T\}$.**

**Thus, if $\#S \geq 2$, then the function $\mu: S \rightarrow S$, such that $\mu(a) = a$ for all $a \in S$, belongs to $C(S) \setminus C^0(S)$.**

Let $\Re$ denote the set of all real numbers and $\Re_+$ the set of non-negative real numbers. Let [X] denote the set of all non-empty subsets of X. Members of [X] are called coalitions.

Given $S \in [X]$, let $C(S) = \{\mu / \mu$ is a bijection on X with $\mu(S) = S\}$ and $C^0(S) = \{\mu \in C(S)$ / T is a non-empty proper subset of S implies $\mu(T) \neq T\}$.

A Contract Choice Problem (CCP) G is an ordered pair $< X, (F(S))_{S \in [X]} >$ such that for all $S \in [X]$: (i) F(S) is a non-empty finite subset of $\Re^S$; (ii) $F(\{a\}) = \{0\}$.

For $G = \langle X, (F(a,b))_{S \in [X]} \rangle$ and $S \in [X]$, $F(S)$ is the set of all feasible allocations of money for agents in S.

A CCP $G = \langle X, (F(S))_{S \in [X]} \rangle$ is said to be super-additive if for all $S, T \in [X]$, with $S \cap T = \phi$: $[x \in F(S), y \in F(T)]$ implies $[z \in F(S \cup T)$, where $z(a) = x(a)$ for all $a \in S$ and $z(a) = y(a)$ for all $a \in T]$.

A CCP $G = \langle X, (F(S))_{S \in [X]} \rangle$ is said to be a generalized matching problem if for all $a \in X$ there exists a function $u^a : X \to \Re$ satisfying the following property: for all $S \in [X]$, $F(S) \subset \{x \in \Re^S /$ for some $\mu \in C^0(S)$, $x(a) = u^a(\mu(a))$ for all $a \in S\} \cup \{-e^S\}$.
The requirement that a generalized matching problem $G = \langle X, (F(S))_{S \in [X]} \rangle$ is a CCP, implies that $F(\{a\}) = \{0\}$ for all $a \in S$. Thus, $F(\{a\}) = \{u^a(a)\}$ implies, $u^a(a) = 0$ for all $a \in S$.

A CCP $G = \langle X, (F(S))_{S \in [X]} \rangle$ is said to be a Shapley-Scarf housing market if for all $a \in X$ there exists a function $u^a : X \to \Re$ satisfying the following property: for all $S \in [X]$, $F(S) = \{x \in \Re^S /$ for some $\mu \in C^0(S)$, $x(a) = u^a(\mu(a))$ for all $a \in S\}$.
Clearly a Shapley-Scarf housing market is a generalized matching problem.

Given a CCP $G = \langle X, (F(S))_{S \in [X]} \rangle$, a coalition structure for G is a partition of X.
A pay-off function is a function $v : X \to \Re_+$. If v is a pay-off function and $S \in [X]$, then $v|S$ denotes the restriction of v to the set S.
An outcome for a CCP $G = \langle X, (F(S))_{S \in [X]} \rangle$ is a pair (f, v), where f is a coalition structure for G and v is a pay-off function such that (i) for all $a \in X$: $v(a) \geq 0$; (ii) for all $S \in f$: $v|S \in F(S)$.
The pair (f, v), where $f = \{\{a\}/ a \in X\}$ and $v(a) = 0$ for all $a \in X$, is an outcome for every CCP. Hence the set of outcomes is always non-empty.

Given $S \in [X]$, let $e^S$ denote the vector in $Z^S$ such that $e^S(i) = 1$ for all $i \in S$ and let #S denote the number of elements of S.
A special case of a CCP is the room-mates problem of Gale and Shapley (1962), where $F(S) = \{-e^S\}$, whenever #S > 2. The marriage problem of Gale and Shapley (1962) is in turn a special case of their room-mates problem. If $F(S) = \{-e^S\}$, whenever #S > 3, then we have a possible generalization of the man, woman and child problem of Alkan (1988).

3. The non-emptiness of the core:

*Given an outcome (f, v) for a CCP $G = \langle X, (F(S))_{S \in [X]} \rangle$, a coalition $S \in [X]$ is said to block (f, v) if there exists $x \in F(S)$: $x(a) > v(a)$ for all $a \in S$.*

*An outcome (f, v) for a CCP $G = \langle X, (F(S))_{S \in [X]} \rangle$ is said to belong to the core of G, if it does not admit any blocking coalition. Let Core(G) denote the set of outcomes in the core of G.*

*An outcome (f, v) for a CCP $G = \langle X, (F(S))_{S \in [X]} \rangle$ is said to be Weakly Pareto Optimal if it does not admit X as a blocking coalition.*

*Given a CCP $G = \langle X, (F(S))_{S \in [X]} \rangle$, an outcome (f, v) is said to be weakly blocked by a coalition $T \subset X$, if there exists $x \in F(T)$: $x(a) \geq v(a)$ for all $a \in T$, with* *strict inequality for at least one $a \in T$. If an outcome (f,v) is weakly blocked by a coalition $T \subset X$, via $x \in F(T)$, then $a \in T$ is said to be an active member of the weakly blocking coalition T, if $x(a) > v(a)$.*

*An outcome (f,v) is said to belong to the strict core of a CCP, if it is not weakly blocked by any coalition. Let SCore(G) denote the set of all outcomes in the strict core of the CCP G. Clearly, Score(G) $\subset$ Core (G).*

*An outcome (f,v) is said to be Pareto Optimal if it does not admit X as a weakly blocking coalition.*

The following result due to Roth and Postelwaite [1977] is well known:
*If G is a Shapley-Scarf housing market, then there exists at least one outcome belonging to the core of G.*

The following example due to Gale and Shapley (1962) shows that the core of a room-mate problem may be empty.

Example 1 (Gale Shapley (1962)) : Let $X = \{1,2,3,4\}$. For $a \in X$, let $u^a: X \to \Re$ be defined as follows:
$u^1$: $u^1(2) = 3$, $u^1(3) = 2$, $u^1(4) = 1$, $u^1(1) = 0$;
$u^2$: $u^2(3) = 3$, $u^2(1) = 2$, $u^2(4) = 1$, $u^2(2) = 0$;
$u^3$: $u^3(1) = 3$, $u^3(2) = 2$, $u^3(4) = 1$, $u^3(3) = 0$;
$u^4$: $u^4(1) = 3$, $u^4(2) = 2$, $u^4(3) = 1$, $u^4(4) = 0$.
Let, $G = <X, (F(S))_{S \in [X]}>$ be the generalized matching problem such that for all $S \in [X]$: (i)$F(S) = \{x \in \Re^S/$ for some $\mu \in C^0(S)$, $x(a) = u^a(\mu(a))$ for all $a \in S\}$, if #S $\in \{1,2\}$; (ii) $F(S) = \{-e^S\}$, otherwise.
Suppose (f,v) is an outcome such that $v(4) \neq 0$. If $v(4) = 1$, then $\{3,4\} \in f$ and $v(3) = 1$. Thus, $\{2,3\}$ blocks (f,v), since 2 can get 3 units and 3 can get 2 units in $F(\{2,3\})$; if $v(4) = 2$, then $\{2,4\} \in f$ and $v(2) = 1$. Thus, $\{1,2\}$ blocks (f,v) since 1 can get 3 units and 2 can get 2 units in $F(\{1,2\})$; if $v(4) = 3$, then $\{1,4\} \in f$ and $v(1) = 1$. Thus, $\{1,3\}$ blocks (f,v) since 3 can get 3 units and 1 can get 2 units in $F(\{1,3\})$. Thus, $v(4) \neq 0$ implies (f,v) does not belong to Core(G). Hence suppose $v(4) = 0$. If $v(3) = 0$, then both $\{2,3\}$ and $\{3,4\}$ block (f,v); if $v(2) = 0$, then both $\{1,2\}$ and $\{2,4\}$ block (f,v); if $v(1) = 0$, then both $\{1,3\}$ and $\{1,4\}$ block (f,v).
Since $v(4) = 0$ requires $v(a) = a$ for at least one $a \in \{1,2,3\}$, Core(G) = $\phi$.
It is worth observing that G is a super-additive CCP.

A sufficient condition for the existence of a non-empty core of a CCP G can be easily obtained, along the lines suggested by Banerjee, Konishi and Sonmez (2001).

Given a CCP G, and a non-empty subset V of X, an outcome (f,v) is said to have the weak top- coalition property for V if there exists a non-empty subset S of V which has a finite partition $\{S^1,\ldots,S^g\}$ satisfying the following properties:
  a. $S \in f$;
  b. For all $a \in S^1$: $v(a) \geq x(a)$ for all $x \in F(T \cup \{a\})$, $T \subset V \setminus \{a\}$;

   c. For all $t \in \{2,\ldots,g\}$, $a \in S^t$, $x \in F(T \cup \{a\})$, $T \subset V \setminus \{a\}$ and $x(a) > v(a)$ implies $T \cap S^k \neq \phi$ for some $k < t$.

*Note: If (f,v) has the weak top coalition property for V with S as defined above, and (f',v') is an outcome, such that $S \in f'$ and $v'(a) = v(a)$ for all $a \in S$, then (f',v') also has the weak top coalition property for V.*

A CCP $G = <X, (F(S))_{S \in [X]}>$ is said to satisfy the weak top coalition property if for any non-empty subset V of X, there exists an outcome (f,v), satisfying the weak top coalition property for V.

Given a CCP $G = <X, (F(S))_{S \in [X]}>$, say that an outcome (f,v) is the union under substitutions of the outcomes $(f^1,v^1),\ldots,(f^g,v^g)$, if: (i) $f = \{S^1,\ldots,S^g\}$; (ii) $S^t \in f^t$, for $t=1,\ldots,g$; (iii) $v|S^t = v^t|S^t$ for $t = 1,\ldots,g$.
If the outcome (f,v) is the union under disjoint substitutions of the outcomes $(f^1,v^1),\ldots,(f^g,v^g)$, then we write $(f,v) = \vee\{(f^1,v^1),\ldots,(f^g,v^g)\}$.

Theorem 1: Let $G = <X, (F(S))_{S \in [X]}>$ be a CCP. If G satisfies the weak top coalition property, then $Core(G) \neq \phi$.

Proof: Let $G = <X, (F(S))_{S \in [X]}>$ be a CCP satisfying the weak top coalition property. Hence, there exists an outcome $(f^1,v^1)$ satisfying the weak top coalition property for X. Let $S_1$ be the non-empty subset of X, having the partition $\{S_1(1),\ldots, S_1(g(1))\}$ such that:
a. $S_1 \in f^1$;
b. For all $a \in S_1(1)$: $v(a) \geq x(a)$ for all $x \in F(T \cup \{a\})$, $T \subset X \setminus \{a\}$;
c. For all $t \in \{2,\ldots,g(1)\}$, $a \in S_1(t)$, $x \in F(T \cup \{a\})$, $T \subset X \setminus \{a\}$ and $x(a) > v(a)$ implies $T \cap S_1(h) \neq \phi$ for some $h < t$.
Hence, no member of $S_1$ will belong to a coalition which blocks $(f^1, v^1)$.
Having defined $S_1,\ldots,S_k$, $(f^1,v^1),\ldots,(f^k, v^k)$ for $k \geq 1$, such that no member of $\bigcup_{j=1}^{k} S_j$ will belong to a coalition which blocks $(f^k, v^k)$, let $(f,v) = (f^k,v^k)$ if $\bigcup_{j=1}^{k} S_j = X$.

If $\bigcup_{j=1}^{k} S_j = X$, then $(f,v) \in Core(G)$. Hence suppose, $\bigcup_{j=1}^{k} S_j \neq X$. Thus, there exists an outcome $(f^{k+1},v^{k+1})$ with $v^{k+1}(a) = v^k(a)$ for all $a \in \bigcup_{j=1}^{k} S_j$ and $\{S_1,\ldots,S_k\} \subset f^{k+1}$, satisfying the weak top coalition property for $X \setminus \bigcup_{j=1}^{k} S_j$. Let $S_{k+1}$ be the non-empty subset of $X \setminus \bigcup_{j=1}^{k} S_j$, having the partition $\{S_{k+1}(1),\ldots, S_{k+1}(g(k+1))\}$ such that:

a. $S_{k+1} \in f^{k+1}$;
b. For all $a \in S_{k+1}(1)$: $v(a) \geq x(a)$ for all $x \in F(T \cup \{a\})$, $T \subset X \setminus \{a\}$;
c. For all $t \in \{2,\ldots,g(k+1)\}$, $a \in S_{k+1}(t)$, $x \in F(T \cup \{a\})$, $T \subset X \setminus \{a\}$ and $x(a) > v(a)$ implies $T \cap S_{k+1}(h) \neq \phi$ for some $h < t$.

Clearly, no member of $\bigcup_{j=1}^{k+1} S_j$ will belong to a coalition which blocks $(f^{k+1}, v^{k+1})$.

Since X is finite, there exists a positive integer K, such that $\bigcup_{j=1}^{K} S_j = X$. Let $(f,v) = (f^K, v^K)$. Thus, $(f,v) \in Core(G)$. Thus, Core(G) is non-empty. Q.E.D.

In fact the following adaptation of yet another property in Banerjee, Konoshi and Sonmez (2001) guarantees that the strict core is a singleton.

Given a CCP $G = <X, (F(S))_{S \in [X]}>$ and a non-empty subset V of X, an outcome $(f,v)$ is said to have the top- coalition property for V if there exists a non-empty subset S of V satisfying the following properties:
  a. $S \in f$;
  b. For all $a \in S$, $T \subset V \setminus \{a\}$ and $x \in F(T \cup \{a\})$: $v(a) \geq x(a)$.

A CCP $G = <X, (F(S))_{S \in [X]}>$ is said to satisfy the top coalition property if for any non-empty subset V of X, there exists an outcome $(f,v)$ satisfying the top coalition property for V.

Theorem 2: Let $G = <X, (F(S))_{S \in [X]}>$ be a CCP. If G satisfies the top coalition property, then $[(f,v), (f',v') \in SCore(G)]$ implies $[v' = v]$.

Proof: Let $G = <X, (F(S))_{S \in [X]}>$ be a CCP satisfying the top coalition property. Hence, there exists an outcome $(f^1, v^1)$ satisfying the top coalition property for X. Let $S_1$ be the non-empty subset of X such that:
a. $S_1 \in f^1$;
b. For all $a \in S_1$: $v(a) \geq x(a)$ for all $x \in F(T \cup \{a\})$, $T \subset X \setminus \{a\}$.

Hence, no member of $S_1$ will be an active member of a coalition which weakly blocks $(f^1, v^1)$.

Having defined $S_1, \ldots, S_k$, $(f^1, v^1), \ldots, (f^k, v^k)$ for $k \geq 1$, such that no member of $\bigcup_{j=1}^{k} S_j$ will be an active member of a coalition which weakly blocks $(f^k, v^k)$, let $(f,v) = (f^k, v^k)$ if $\bigcup_{j=1}^{k} S_j = X$. If $\bigcup_{j=1}^{k} S_j = X$, then $(f,v) \in SCore(G)$. Hence suppose, $\bigcup_{j=1}^{k} S_j \neq X$. Thus, there exists an outcome $(f^{k+1}, v^{k+1})$ with $v^{k+1}(a) = v^k(a)$ for all $a \in \bigcup_{j=1}^{k} S_j$ and $\{S_1, \ldots, S_k\} \subset f^{k+1}$, satisfying the top coalition property for $X \setminus \bigcup_{j=1}^{k} S_j$. Let $S_{k+1}$ be the non-empty subset of $X \setminus \bigcup_{j=1}^{k} S_j$, such that:

a. $S_{k+1} \in f^{k+1}$;
b. For all $a \in S_{k+1}$: $v(a) \geq x(a)$ for all $x \in F(T \cup \{a\})$, $T \subset X \setminus \{a\}$;

Clearly, no member of $\bigcup_{j=1}^{k+1} S_j$ will be an active member of a coalition which weakly blocks $(f^{k+1}, v^{k+1})$.

Since X is finite, there exists a positive integer K, such that $\bigcup_{j=1}^{K} S_j = X$. Let $(f,v) = (f^K, v^K)$. Clearly no member of X can be an active member of a coalition which weakly blocks (f,v). Thus, $(f,v) \in SCore(G)$.

Let $(f',v') \in SCore(G)$. If $v'(a) \neq v(a)$ for some $a \in S^1$, then $S^1$ can weakly block $(f',v')$. Thus, $(f',v') \in SCore(G)$ implies $v'(a) = v(a)$ for all $a \in S^1$.

Suppose K > 1. Suppose, $v'(a) = v(a)$ for all $a \in \bigcup_{t=1}^{k} S^t$, for some $k \in \{1,\ldots,K-1\}$. If $v'(a) \neq v(a)$ for some $a \in S^{k+1}$, then $S^{k+1}$ can weakly block $(f',v')$. Thus, $(f',v') \in SCore(G)$ implies $v'(a) = v(a)$ for all $a \in \bigcup_{t=1}^{k+1} S^t$. Hence, $(f',v') \in SCore(G)$ implies $v' = v$. Q.E.D.

4. Existence of Individually Stable Outcomes:

*Given a CCP $G = < X, (F(S))_{S \in [X]} >$ an outcome (f,v) is said to be unilaterally blocked by agent $a \in X$, if there exists a non-empty subsets S, T of X with $S\setminus\{a\} \subset T$, $T \in f$ and $x \in F(S)$ such that $x(b) > v(b)$ for all $b \in S$.*

*An outcome (f,v) is said to be individually stable for G if it is not blocked by any agent. Let IS(G) denote the set of all individually stable outcomes for G.*

Given a non-empty subset V of X, an outcome (f,v) is said to have the weak top-cycle property for V if there exists a non-empty subset S of V containing 's' distinct elements and a one-to-one function $\psi: S \to \{1,\ldots,s\}$ satisfying the following properties:

1. $S \in f$;
2. For all $a \in S$, for all non-empty subsets T of $V \setminus S$ and $x \in F(T \cup \{a\})$: $v(a) \geq x(a)$;
3. For all $a,b \in S$ with $\psi(a) < \psi(b)$, for all non-empty subsets T of $V \setminus \{a\}$ with $b \in T$ and $x \in F(T \cup \{a\})$: $[x(a) > v(a)]$ implies $[v(c) \geq x(c)$ for some $c \in T]$.

A CCP $G = < X, (F(S))_{S \in [X]} >$ is said to satisfy the weak top cycle property if for any non-empty subset V of X, there exists an allocation (f,v), satisfying the weak top cycle property for V.

Theorem 3: Let G be a CCP satisfying the weak top cycle property, then IS(G) $\neq \phi$.

Proof: Let $G = < X, (F(S))_{S \in [X]} >, (rk_i)_{i \in N} >$ be a CCP satisfying the weak top cycle property. Thus, there exists an outcome $(f^1,v^1)$ satisfying the weak top cycle property for X. Let $S^1$ be a non-empty subset of X, containing $s^1$ distinct agents and a bijection $\psi^1: S^1 \to \{1,\ldots,s^1\}$ such that: (a) $S^1 \in f$; (b) For all $a \in S^1$, for all non-empty subsets T of $X \setminus S^1$ and $x \in F(T \cup \{a\})$: $v(a) \geq x(a)$; (c) For all $a,b \in S^1$ with $\psi^1(a) < \psi^1(b)$, for all non-empty subsets T of $V \setminus \{a\}$ with $b \in T$ and $x \in F(T \cup \{a\})$: $[x(a) > v(a)]$ implies $[v(c) \geq x(c)$ for some $c \in T]$.

Having obtained $(f^t, v^t)$, $S^t$, $\psi^t$ for $\tau \geq t \geq 1$, if $X \setminus \bigcup_{t=1}^{\tau} S^t \neq \phi$, let $(f^{\tau+1}, v^{\tau+1})$ an outcome, $S^{\tau+1}$ a non-empty subset of $X \setminus \bigcup_{t=1}^{\tau} S^t$ containing $s^{\tau+1}$ distinct agents and a bijection $\psi^{\tau+1}: S^{\tau+1} \to \{1,\ldots,s^{\tau+1}\}$ be such that: (a) $S^{\tau+1} \in f^{\tau+1}$; (b) For all $a \in S^{\tau+1}$, for all non-empty subsets T of $X \setminus \bigcup_{t=1}^{\tau+1} S^t$ and $x \in F(T \cup \{a\})$: $v(a) \geq x(a)$; (c) For all $a,b \in S^{\tau+1}$ with $\psi^1(a) < \psi^1(b)$, for all non-empty subsets T of $V \setminus \{a\}$ with $b \in T$ and $x \in F(T \cup \{a\})$: $[x(a) > v(a)]$ implies $[v(c) \geq x(c)$ for some $c \in T]$.

Since X is a non-empty finite set, there exists a least positive integer g, such that $\bigcup_{t=1}^{g} S^t = X$. Let (f,v) be the outcome such that $f = \{S^1,\ldots, S^g\}$ and for all $k \in \{1,\ldots,g\}$: $v(a) = v^k(a)$, for all $a \in S^k$.

Suppose there is a coalition S, such that $S \setminus \{a\} \subset T$ for some $T \in f$ and $x \in F(S)$ such that: $x(c) > v(c)$ for all $c \in S$. Suppose $a \in S^1$. Since, $a \in S^1$ implies $v(a) \geq x(a)$ if $S \setminus \{a\} \subset X \setminus S^1$, it must be the case that $S \setminus \{a\} \cap S^1 \neq \phi$. If there exists $b \in S \setminus \{a\} \cap S^1$, such that $\psi^1(b) > \psi^1(a)$, then $v(c) \geq x(c)$ for some $c \in S \setminus \{a\}$. Hence suppose, $\psi^1(a) > \psi^1(b)$ for all $b \in S \setminus \{a\} \cap S^1$. Let $b^*$ be the unique element in S such that $\psi^1(b^*) \leq \psi^1(c)$ for all $c \in S \cap S^1$. Thus, $x(b^*) > v(b^*)$, $\psi^1(b^*) < \psi^1(a)$, $a,b \in S^1$, implies that there exists $c \in S \setminus \{b^*\}$ such that $v(c) \geq x(c)$. This contradicts $x(c) > v(c)$ for all $c \in S$.

Now suppose, $a \in X \setminus S^1$ and $S \cap S^1 \neq \phi$. Let $b^*$ be the unique element in S such that $\psi^1(b^*) \leq \psi^1(c)$ for all $c \in S \cap S^1$. Suppose, $S \setminus \{b^*\} \cap S^1 \neq \phi$. Let $b \in S \setminus \{b^*\} \cap S^1$. Thus, $x(b^*) > v(b^*)$, $\psi^1(b^*) < \psi^1(b)$, $b^*,b \in S^1$, implies that there exists $c \in S \setminus \{b^*\}$ such that $v(c) \geq x(c)$. This contradicts $x(c) > v(c)$ for all $c \in S$. Thus, let $S \setminus \{b^*\} \cap S^1 = \phi$. Then, $v(b^*) \geq x(b^*)$ leads to a contradiction once again.

Hence, no agent in $S^1$ will unilaterally deviate from (f,v), and any agent who unilaterally deviates from (f,v), must belong to $X \setminus S^1$ and the coalition that he would be joining to deviate must also be a subset of $X \setminus S^1$.

Proceeding as above, it is easily observed that for all $k \in \{1,\ldots,g\}$, no a belonging to $\bigcup_{t=0}^{k} S^t$ will unilaterally deviate from (f,v). Thus, $(f,v) \in IS(G)$. Hence, $IS(G) \neq \phi$.
Q.E.D.

6.     Non-emptiness of the weak bargaining set:

*Given a CCP $G = \langle X, F(S)_{S \in [X]} \rangle$ and a Pareto Optimal outcome $(f,v)$, the pair $((f',v'), T)$ where $(f',v')$ is an outcome for G and $T \in f'$ is said to be a strong objection against $(f,v)$ if $v'(a) > v(a)$ for all $a \in T$ and no subset of T is a blocking coalition for $(f',v')$.*

*Given a CCP $G = \langle X, F(S)_{S \in [X]} \rangle$, a Pareto Optimal outcome $(f,v) \in F$ and a strong objection $(f',v'),T$ against $(f,v)$, an ordered pair $((f'',v''),U)$ where $(f'',v'')$ is an outcome for G and $U \in f''$ is said to be a strong counter-objection against $((f',v'),T)$ if: (a) $U \setminus T$, $T \setminus U$ and $U \cap T$ are all non-empty; (b) $v''(a) > v'(a)$ for all $a \in U$.*

*The strong objection $((f',v'),T)$ against the outcome $(f,v)$ is said to be justified, if $((f',v'),T)$ has no strong counter-objection.*

*We define the weak bargaining set of a CCP $G = \langle X, F(S)_{S \in [X]} \rangle$, to be the set $WB(G) = \{(f,v) / (f,v)$ is Pareto Optimal and such that no strong objection against $(f,v)$ is justified$\}$.*

Example 1 (due to Gale and Shapley (1962)) is one which has an empty core, but a non-empty weak bargaining set. As in Example 1, let $X = \{1,2,3,4\}$. For $a \in X$, let $u^a: X \to \Re$ be defined as follows:
$u^1$: $u^1(2) = 3$, $u^1(3) = 2$, $u^1(4) = 1$, $u^1(1) = 0$;
$u^2$: $u^2(3) = 3$, $u^2(1) = 2$, $u^2(4) = 1$, $u^2(2) = 0$;
$u^3$: $u^3(1) = 3$, $u^3(2) = 2$, $u^3(4) = 1$, $u^3(3) = 0$;
$u^4$: $u^4(1) = 3$, $u^4(2) = 2$, $u^4(3) = 1$, $u^4(4) = 0$.
Let, $G = \langle X, (F(S))_{S \in [X]} \rangle$ be the generalized matching problem such that for all $S \in [X]$: (i) $F(S) = \{x \in \Re^S /$ for some $\mu \in C^0(S)$, $x(a) = u^a(\mu(a))$ for all $a \in S\}$, if $\#S \in \{1,2\}$; (ii) $F(S) = \{-e^S\}$, otherwise.
We saw in Example 1, that $Core(G) = \phi$.
However, consider $v(4) = 1$, $v(3) = 1$, $v(2) = 2$, $v(1) = 3$, $f = \{\{1,2\},\{3,4\}\}$. The pair $((f', v'), \{2,3\})$ is a strong objection against $(f,v)$, where $v'(2) = 3$, $v'(3) = 2$, $v'(1) = v'(4) = 0$ and $f' = \{\{1\},\{4\},\{2,3\}\}$. Let $f'' = \{\{2\}, \{4\}, \{1,3\}\}$, $v''(1) = 2$, $v''(3) = 3$, $v''(2) = v''(4) = 0$. Then the pair $((f'',v''), \{1,3\})$ is a strong counter-objection against $((f',v'),\{2,3\})$. Further, $(f,v)$ admits no blocking coalition other than $\{1,3\}$. Since no subset of $\{1,3\}$ blocks $(f',v')$, $(f,v)$ belongs to $WB(G)$.
Note that the outcome $(f^*, v^*)$ such that $f^* = \{\{1,2,3\},\{4\}\}$ and $v^*(1) = v^*(2) = v^*(3) = 3$, $v^*(4) = 0$, belongs to the $Core(G^*)$, where $G^* = \langle X, F(S)_{S \in [X]} \rangle$ is such that for all $S \in [X]$: (i) $F(S) = \{x \in \Re^S /$ for some $\mu \in C^0(S)$, $x(a) = u^a(\mu(a))$ for all $a \in S\}$, if $\#S \in \{1,2,3\}$; (ii) $F(S) = \{-e^S\}$, otherwise.

Note that it is possible to provide a definition of the weak bargaining set modified along the lines suggested in Mas-Colell (1989).

*Given a CCP $G = \langle X, F(S)_{S \in [X]} \rangle$ an outcome $(f,v)$ and a strong objection $((f',v'),T)$ against $(f,v)$, an ordered pair $((f'',v''),U)$ is said to be a classical strong counter-objection against $((f',v'),T)$ if: (a) $U \in f''$; (b) $U \setminus T$, $U \setminus S$ and $U \cap T$ are all non-empty; (c) $v''(a) \geq v(a)$ for all $a \in U \setminus T$; (d) $v''(a) > v'(a)$ for all $a \in U$.*

*The strong objection ((f',v'),T) against the outcome (f,v) is said to be classically justified, if ((f',v'),T) has no classical strong counter-objection.*
*We define the classical weak bargaining set of a CCP G = <X,F(S)$_{S\in[X]}$>, to be the set WB$^*$(G) = {(f,v)/ (f,v) is Pareto Optimal, and such that no strong objection against (f,v) is classically justified}.*

However, the following example reveals that even for room-mates problems, WB$^*$(G) may be empty.

Example 2: Let X = {1,2,3}. For a∈X, let $u^a$: X→$\Re$ be defined as follows:
$u^1$: $u^1(2) = 2$, $u^1(3) = 1$, $u^1(1) = 0$;
$u^2$: $u^2(3) = 2$, $u^2(1) = 1$, $u^2(2) = 0$;
$u^3$: $u^3(1) = 3$, $u^3(2) = 2$, $u^3(3) = 0$.

Let, G = < X, (F(S))$_{S\in[X]}$ > be the generalized matching problem such that for all S∈[X]: (i)F(S) = {x∈$\Re^S$/ for some µ∈$C^0$(S), x(a) = $u^a$(µ(a)) for all a∈S}, if #S ∈{1,2}; (ii) F(S) ={-$e^S$}, otherwise.
Since (f,v) such that v(a) = 0 for all a∈X is not Pareto Optimal, it cannot belong to WB$^*$(G).
Let (f,v) be the outcome such that f = {{1,3},{2}}, v(1) = 1, v(2) = 0, v(3) = 2 and (f',v') be the outcome such that f' = {{1,2},{3}}, v'(1) = 2, v'(2) = 1, v'(3) = 0. Thus, ((f',v'), {1,2}) is a strong objection against (f,v). Any strong counter-objection or classical strong counter-objection cannot contain agent 1, since agent 1 gets 2 units of money at (f',v'). The only possibility is (({{2,3}, {1}}, v"), {2,3}) where v"(1) = 0, v"(3) = 1, v"(2) = 2, which is a strong counter-objection though not a classical strong counter-objection, since agent 3 is worse off at (f",v") than at (f,v). Thus, (f,v) ∉WB$^*$(G).
Let (f,v) be the outcome such that f = {{1},{2,3}}, v(1) = 0, v(2) = 2, v(3) = 1 and (f',v') be the outcome such that f' = {{1,3},{2}}, v'(1) = 1, v'(2) = 0, v'(3) = 2. Thus, ((f',v'), {1,3}) is a strong objection against (f,v). Any strong counter-objection or classical strong counter-objection cannot contain agent 3, since agent 3 gets 2 units of money at (f',v'). The only possibility is (({{1,2}, {3}}, v"), {1,3}) where v"(1) = 2, v"(3) = 0, v"(2) = 1, which is a strong counter-objection though not a classical strong counter-objection, since agent 2 is worse off at (f",v") than at (f,v). Thus, (f,v) ∉WB$^*$(G).
Hence Bar$^*$(G) = φ.

Theorem 4: Let G be a CCP. Then, WB(G) ≠ φ.

Proof: Let G = < X, F(S)$_{S\in[X]}$> be a CCP and let (f,v) be a Pareto Optimal outcome for G. If (f,v) does not admit a strong objection then clearly, (f,v) ∈WB(G). Suppose (($f^1,v^1$), $S^1$) is a strong objection against (f,v) which further does not admit a strong counter-objection. Then, no member of $S^1$ is part of a strong objection against ($f^1,v^1$). Clearly, there can be no strong objection (($f^2,v^2$), $S^2$) against ($f^1,v^1$) such that $S^2 \cap S^1 \neq \phi$. If ($f^1,v^1$) does not admit any strong objection, then ($f^1,v^1$)∈WB(G). Suppose (($f^2,v^2$), $S^2$) is a strong objection against

($f^1,v^1$) which further does not admit a strong counter-objection. Clearly, $S^2 \cap S^1 = \phi$. Without loss of generality suppose $S^1 \in f^2$ and $v^2(a) = v^1(a)$ for all $a \in S^1$. This is possible, since $S^2 \cap S^1 = \phi$. Then, no member of $S^1 \cup S^2$ is part of a strong objection against ($f^2,v^2$).

Having constructed a strong objections (($f^p,v^p$), $S^p$) against ($f^{p-1},v^{p-1}$) for $p = 1,\ldots,$ k, where ($f^0,v^0$) = (f,v), such that no member of $\bigcup_{p=1}^{k} S^p$ is part of a blocking coalition against ($f^k,v^k$) there are two possibilities: there does exist a strong objection against ($f^k,v^k$) in which case ($f^k,v^k$)$\in$WB(G); there exists a strong objection (($f^{k+1}, v^{k+1}$), $S^{k+1}$) against ($f^k,v^k$). If every such strong objection admits a strong counter-objection, then ($f^k,v^k$)$\in$WB(G). If not then there exists a strong objection (($f^{k+1},v^{k+1}$), $S^{k+1}$), which further does not admit a strong counter-objection. Clearly, $S^{k+1} \cap (\bigcup_{p=1}^{k} S^p) = \phi$. Without loss of generality suppose, $S^p \in f^{k+1}$ for $p = 1,\ldots,k$ and $v^{k+1}(a) = v^k(i)$ for all $a \in \bigcup_{p=1}^{k} S^p$. Then no member of $\bigcup_{p=1}^{k+1} S^p$ is part of a strong objection against ($f^{k+1}, v^{k+1}$).

Since X is a finite set, there is a smallest positive integer K, such that either every objection ((f',v'), T) against ($f^K,v^K$) admits a strong counter-objection, or [$\bigcup_{p=1}^{K} S^p$ = X or no member of X is part of a blocking coalition against ($f^K,v^K$). In either case, ($f^K,v^K$)$\in$WB(G). Q.E.D.